\documentclass[a4paper]{amsart} 

\usepackage{amsmath,amsthm,amssymb,amsfonts,mathrsfs,color,hyperref, mathtools,crop, graphicx, enumitem}
\usepackage{todonotes}
\usepackage{color}

\theoremstyle{plain}
\begingroup
\newtheorem{theorem}{Theorem}[section]

\endgroup

\theoremstyle{definition}
\begingroup

\endgroup

\theoremstyle{remark}
\begingroup
\newtheorem{remark}[theorem]{Remark}
\endgroup 

\numberwithin{equation}{section}
\newcommand{\N}{\mathbb N} 
 
\newcommand{\R}{\mathbb R} 

\newcommand{\dist}{{\rm dist}}
\newcommand{\diam}{{\rm diam}}

\newcommand{\wto}{\rightharpoonup}

\newcommand{\W}{{\mathcal W}}

\renewcommand{\L}{{\mathcal L}}

\newcommand{\C}{{\mathcal C}}
\newcommand{\F}{{\mathcal F}}

\newcommand{\LRa} {\Leftrightarrow}
\newcommand{\Ra} {\Rightarrow}

\renewcommand{\H}{{\mathcal H}}
\newcommand{\spt}{{\operatorname{spt}}}
\newcommand{\inv}{^{-1}}
\newcommand{\cc}{\Subset}

\newcommand{\eps}{\varepsilon}

\newcommand{\argmin}{\mathrm{argmin}}

\newcommand{\dx}{\,\mathrm{d}x}
\renewcommand{\d}{\,\mathrm{d}}
\newcommand{\dy}{\,\mathrm{d}y}
\newcommand{\ds}{\,\mathrm{d}s}

\def\Xint#1{\mathchoice
{\XXint\displaystyle\textstyle{#1}}%
{\XXint\textstyle\scriptstyle{#1}}%
{\XXint\scriptstyle\scriptscriptstyle{#1}}%
{\XXint\scriptscriptstyle\scriptscriptstyle{#1}}%
\!\int}
\def\XXint#1#2#3{{\setbox0=\hbox{$#1{#2#3}{\int}$ }
\vcenter{\hbox{$#2#3$ }}\kern-.6\wd0}}

\def\dashint{\Xint-}

\begin{document}

\title[Keeping it together]{Keeping it together: a phase field version of path-connectedness and its implementation}

\author{Patrick W.~Dondl}
\address{Patrick W.~Dondl\\
Abteilung f\"ur Angewandte Mathematik\\
Albert-Ludwigs-Universit\"at Freiburg\\
Hermann-Herder-Str.~10\\
79104 Freiburg i.~Br.\\
Germany}
\email{patrick.dondl@mathematik.uni-freiburg.de}

\author{Stephan Wojtowytsch}
\address{Stephan Wojtowytsch\\Department of Mathematical Sciences\\ Carnegie Mellon University\\ 5000 Forbes Avenue\\ Pittsburgh, PA 15213\\ USA}
\email{swojtowy@andrew.cmu.edu}

\date{\today}

\subjclass[2010]{49M30, 90C59}
\keywords{Willmore Energy, Phase Field Approximation, Topological Constraint, Connectedness, Dijkstra's Algorithm}

\begin{abstract}
We describe the implementation of a topological constraint in finite element simulations of phase field models which ensures path-connectedness of preimages of intervals in the phase field variable. Two main applications of our method are presented. First, a discrete steepest decent of a phase field version of a bending energy with spontaneous curvature and additional surface area penalty is shown, which leads to disconnected surfaces without our topological constraint but connected surfaces with the constraint. The second application is the segmentation of an image into a connected component and its exterior. Numerically, our constraint is treated using a suitable geodesic distance function which is computed using Dijkstra's algorithm.
\end{abstract}

\maketitle

\section{Introduction}
In this article we describe how to incorporate a topological constraint into a phase field simulation for certain geometric functionals. In three dimensions, the prototypical example of such an energy is Willmore's energy, i.e., the integral of mean curvature squared, restricted to the class of $C^2$-manifolds embedded into a bounded domain $\Omega\subset \R^3$ such that the embedding has surface area $S>0$ and is \emph{connected}. A number of more general functionals are also admissible---the precise theoretical setting for our methods is presented in section \ref{sec:theory}. We furthermore consider the case of functionals controlling only the perimeter of sets contained in a bounded domain $\Omega\subset\R^2$. To simplify the presentation, we mostly focus on the three-dimensional case with a control of Willmore's energy.

Our method to enforce this connectedness constraint for diffuse interfaces is based on a functional $\C_\eps$ introduced in~\cite{MR3590663} and given by
\[
\C_\eps(u) = \int_\Omega\int_\Omega \frac{\widetilde W(u(x))}\eps\:\frac{\widetilde W(u(y))}\eps\:d^{F(u)}(x,y)\dx\dy
\]
where $\widetilde W, F$ are continuous functions such that
\[
\widetilde W, F\geq 0, \qquad F(z)=0\LRa z\in [\alpha,\beta], \qquad \widetilde W(z) >0\LRa z\in (\alpha,\beta)
\]
for some $-1<\alpha<\beta<1$. For a heuristic interpretation of the functional, see section~\ref{sec:top_term}. 

Despite the intimidating appearance of the functional $\C_\eps$ as a double integral coupled to a geodesic distance, an efficient algorithm to treat this term can be implemented. We rely on a decomposition into connected components and a Dijkstra or fast marching-type method to approximate the geodesic distance as well as its variation.

In our simulations, we compute a discrete time gradient flow (approximately) minimizing in each time step the functional
\begin{equation} \label{eq:time_step}
u_{n+1} \in \argmin \left( u\mapsto \frac{\eps}{2\tau} ||u-u_n||_{L^2(\Omega)}^2 + \F_\eps(u)+ \frac{a}\eps\,\C_\eps(u)\right)
\end{equation}
where $\F_\eps$ is an admissible phase field approximation of a geometric functional $\F$ and $a\ge 0$ is a constant. Note that the case $a=0$ of course provides no penalty for disconnectedness. 

The article is structured as follows. In section \ref{sec:theory}, we recall precise statements regarding the sharp interface limit of functionals of the form $\F_\eps + \frac{1}{\eps}C_\eps$, where $\F_\eps$ is a term controlling a diffuse version of Willmore's energy as well as the perimeter~\cite{bellettini:2009ui, roger:2006ta}. In the remainder of section \ref{sec:preliminaria}, we give a heuristic explanation of our topological functional and briefly review the background material on the distance function on graphs which will be required below. In section \ref{sec:algorithm}, we describe an implementation of the connectedness constraint. In section \ref{sec:numerics}, we present simulations with and without the topological constraint. Section \ref{sec:2d_conn} provides an outlook to the aforementioned application to image segmentation.

\section{Preliminaries}\label{sec:preliminaria}

\subsection{Geometric Energies and Phase-field Connectedness}\label{sec:theory}
We note the following $\Gamma$-convergence results.

\begin{theorem}\cite{MR3590663} \label{thm:conn}
Let $n=2,3$ and $\F_\eps$ a sequence of functionals such that
 for every $C>0$ there exists $C'>0$ independent of $\eps$ such that
\[
\F_\eps(u)\leq C\qquad\Ra\qquad (\W_\eps + S_\eps)(u) \leq C'
\]
where
\[
S_\eps(u) = \frac1{c_0}\int_\Omega\frac\eps2\,|\nabla u|^2 + \frac1\eps\,W(u)\dx, \qquad c_0 = \int_{-1}^1\sqrt{2W(s)\,}\d s
\]
is the usual Modica-Mortola approximation of the perimeter functional, $W(u) = \frac14\,(u^2-1)^2$ and
\[
\W_\eps(u) = \frac1{c_0\,\eps}\int_\Omega \left(-\eps\,\Delta u + \frac1\eps\,W'(u)\right)^2\dx, \qquad u\in -1 + W_0^{2,2}(\Omega)
\]
is an approximation of Willmore's energy due to de Giorgi. If $\mu$ is a Radon measure such that the diffuse area measures
\[
\mu_\eps^{u_\eps} = \frac1{c_0}\left(\frac\eps2\,|\nabla u_\eps|^2 + \frac1\eps\,W(u_\eps)\right)\L^n
\]
converge to $\mu$ in the weak* sense for a sequence $u_\eps$ with $\limsup_{\eps\to 0}\left[\F_\eps+\eps^{-\kappa}\C_\eps\right](u_\eps)<\infty$, then $\spt(\mu)$ is \emph{connected}.
\end{theorem}

\begin{remark}
According to \cite{roger:2006ta}, we have $\Gamma-\lim_{\eps\to 0} \W_\eps = \W$ at $C^2$-boundaries for all sequences such that $S_\eps(u_\eps)$ also remains bounded. Furthermore, if $(\W_\eps + S_\eps)(u_\eps)$ is uniformly bounded, then up to a subsequence we know that there exist $u\in BV(\Omega, \{-1,1\})$ and a Radon measure $\mu$ supported in $\overline\Omega$ such that $u_\eps\to u$ in $L^1(\Omega)$ and $\mu_\eps^{u_\eps} \wto \mu$ weakly in the sense of Radon measures. Thus our assumption that a limit exists is not detrimental to the generality of the statement.
\end{remark}
\begin{remark}
We note that furthermore, if 
\[
\Gamma(L^1)-\lim_{\eps\to 0} \F_\eps = \F
\]
at $C^2$-boundaries and the recovery sequence for $\F_\eps$ is given by the usual optimal profile construction, then for any $\kappa>0$, the functional $\F_\eps + \eps^{-\kappa}\C_\eps$ satisfies
\[
\Gamma(L^1)-\lim_{\eps\to 0} \left[\F_\eps + \eps^{-\kappa} \C_\eps\right] = \F
\]
at {\em connected} $C^2$-boundaries. The $\liminf$-inequality here is obvious since $\F_\eps + \eps^{-\kappa}\C_\eps\geq \F_\eps$. If $\partial E$ is a connected $C^2$-boundary with associated optimal profile sequence $u_\eps$, then $\C_\eps(u_\eps) \equiv 0$ as demonstrated in \cite{MR3590663} and the $\limsup$-inequality follows as well. This shows that our topological functional does not change existing $\Gamma$-limits, apart from enforcing connectedness.
\end{remark}

Our main numerical example for this result is the treatment of the functional
\begin{equation} \label{eq:will_h0}
\F_\eps(u)  = \frac1{c_0\eps} \int_\Omega \left(-\eps\,\Delta u - \frac1\eps\,W'(u) - H_0\,\sqrt{2\,W(u)}\right)^2\dx + \lambda\, S_\eps(u)
\end{equation}
with $\lambda>0$ and $H_0\in\R$ in three ambient space dimensions. Like above, the normalizing constant is given by $c_0=\int_{-1}^{1} \sqrt{2W(s)}\ds = \frac{4\sqrt{2}}{3}$. A modification of the diffuse Willmore functional $\W_\eps$, this energy includes a spontaneous curvature $H_0$ which the surface would prefer to take. This is a more realistic model for many biological membranes. Theorem~\ref{thm:conn} clearly applies to this $\F_\eps$, but the hypotheses are also satisfied by an approximation of Helfrich's Energy as given in \cite{bellettini:2009ui} or if a volume constraint is included, see \cite[Chapter 5]{MR3590663}. We do note, however, that the $\Gamma$-limit of this functional $\F_\eps$ is not simply given by Willmore's energy with spontaneous curvature and area penalty, as that is not a lower-semicontinuous energy~\cite{grosse1993new}.
 
\subsection{The Topological Term} \label{sec:top_term}
When $(\W_\eps+ S_\eps)(u)<\infty$, $u\in W^{2,2}(\Omega)$ has a continuous representative, so the following notions are well-defined. We can think of the interface as the pre-image of any interval $(\alpha,\beta)\cc (-1,1)$, for a precise statement see~\cite[Theorem 2.20]{DW_conv}. Thus it is possible to introduce a quantitative notion of path-connectedness of the interface $I:= u\inv(\alpha,\beta)$ at two points $x,y\in I$ through a geodesic distance function
\[
d^{F(u)}(x,y) = \inf\left\{\int_\gamma F(u)\d\H^1\:\bigg|\: \gamma\text{ curve from }x\text{ to }y\right\}
\]
with a weight $F(u)$ satisfying 
\[
F\in C^1(\R), \qquad F\equiv 0\text{ on }[\alpha,\beta], \qquad F>0\text{ outside }[\alpha,\beta].
\]
In particular, if $I$ is (path-)connected, then $d^{F(u)}(x,y) = 0$ for all $x,y\in I$. If, however, there are multiple connected components of $I$ with a positive spatial separation, then any curve $\gamma$ should have uniformly positive length between $x$ and $y$ if the two points are in different connected components. We now measure the total disconnectedness of $I$ by a double-integral of the quantitative path-disconnectedness of $I$ at $x$ and $y$ over the entire interface with respect to both $x$ and $y$
\[
\C_\eps(u) = \iint_{\Omega\times\Omega} \frac1\eps \widetilde{W}(u(x)) \,\frac1\eps\widetilde{W}(u(y))\,d^{F(u)}(x,y)\dx\dy
\]
where $\widetilde{W}$ is a bump function 
\[
\widetilde{W}\in C^1(\R), \qquad \widetilde{W}>0 \text{ on }(\alpha,\beta), \qquad \widetilde{W}\equiv 0\text{ outside }(\alpha,\beta).
\]
The dependence of $\C_\eps$ on the choices of $F$ and $\widetilde{W}$, and therefore on $\alpha$ and $\beta$, should be kept in mind, but for notational convenience we will not make it explicit in the remainder of this article. Along the usual optimal profile recovery sequence for a connected, smooth, embedded manifold, $\C_\eps$ vanishes identically. The bound on $\W_\eps + S_\eps$ enforces a strong mode of convergence for $u_\eps$ to $\pm 1$ away from $\operatorname{spt} \mu$  which suffices for $\C_\eps$ to detect a disconnected interface in the sense that $\liminf \C_\eps(u_\eps) >0$ if $\mu$ has more than one connected component, see \cite{MR3590663, DW_conv}. The normalization factor $\frac1\eps$ is used since an interface has width proportional to $\eps$.

\subsection{Discretizing the Geodesic Distance} \label{sec:disc_dist}

Let $\Gamma$ be a finite connected (undirected) graph with vertices $v$ and edges $e$ that are assigned weights $w_e\geq 0$. The distance of two vertices $v, v'$ is defined by the length of the shortest path in the graph connecting $v$ and $v'$ where the length of an edge $e$ is measured by the weight $w_e$, i.e.,
\[
d(v, v') = \inf\left\{\sum_{i=1}^nw_{e_i} \:|\: v= v_0, v' = v_n, \quad v_{i-1}, v_i \in e_i, \quad n\in\N\right\}.
\]

In our finite element setting we consider a sequence of quasi-uniform triangulations (or tetrahedralizations in 3D) $\mathcal{T}_h$ of our domain $\Omega$ with a spatial grid scale $h$ for $h\to 0$. We furthermore assume that $u$ is a given continuous function on $\Omega$. Let now $\Gamma_h(u)$ be the dual graph associated to a triangulation $\mathcal{T}_h$, i.e., a vertex $v$ of $\Gamma_h(u)$ corresponds bijectively to a triangle $T = T_v \in \mathcal{T}_h$ and that if $T_1\in\mathcal{T}_h$ shares an edge (or a face in 3D) with $T_2\in\mathcal{T}_h$, then the corresponding vertices in $\Gamma_h(u)$ are connected by an edge. To such an edge $e=e(T_1,T_2)$ we associate the weight $w_e(u_{T_1},u_{T_2}) = \frac{F(u_{T_1}) + F(u_{T_2})}{2}\cdot\frac{\diam(T_1) + \diam(T_2)}{2}$, where $u_T = \frac1{|T|}\int_T u\dx =: \dashint_T u \dx$. This gives rise to a discrete distance $d^{\Gamma_h(u)}$ depending on a phase field function $u$.

The triangulations may force a minimal path to zig-zag to connect two points, so the distance on the graph may not approximate the distance function
\[
d^{F(u)}(x,y) = \inf\left\{\int_\gamma F(u)\d\H^1\:\bigg|\: \gamma\text{ curve from }x\text{ to }y\right\},
\]
as $h\to 0$, but, due to the non-degeneracy assumptions on the triangulation, we note that the discrete distance $d^{\Gamma_h(u)}$ is equivalent to $d^{F(u)}$ in an almost bi-Lipschitz sense uniformly in $h$, i.e.,
\[
c\,d^{\Gamma_h}(T, T')-\bar{c}h \leq d^{F(u)}(x, x') \leq C\,d^{\Gamma_h}(T, T')+\bar{C}h
\]
for all $x, x' \in \Omega$ and $T,T'\in \mathcal{T}_h$ such that $x\in T, x' \in T'$ with suitable constants $c,\bar{c},C,\bar{C}$. Such a modification clearly does not change the effect of the topological term on connectedness as only a non-zero lower bound and a vanishing upper bound are required for the $\liminf$-inequality and $\limsup$-construction, respectively.

As a treatment of the time-step minimization problem will require a variation of the distance as well, so we note that if there exists a unique shortest curve $\bar\gamma$ between $x$ and $y$ then
\[ 
\frac{\d}{\d t}\bigg|_{t=0}d^{F(u+t\phi)}(x,y) = \int_{\bar \gamma}F'(u)\,\phi\d\H^1.
\]
This identity will be postulated for the procedure below and we approximate the variation of the geodesic distance on the graph as we vary $u$ in direction $\phi$ by the discrete term
\begin{align}\label{eq variation of distance}
\delta_{u;\phi}d^{\Gamma_h}(T,T') &= 
 \sum_{j=1}^n \Bigg[ \frac{\diam(T_j) + \diam(T_{j-1})}2 \cdot \\
&\quad\quad\quad\quad \frac{1}{2} \bigg( F'(u_{T_j})\,\dashint_{T_j}\phi\dx+ F'(u_{T_{j-1}})\, \dashint_{T_{j-1}}\phi\dx\bigg)\Bigg] \nonumber
\end{align}
where $\{T_j\}_{j=0}^n$ is a shortest connecting path between triangles $T=T_0$ and $T'=T_n$ containing $x$ and $y$, respectively. The above variation is in the spirit of~\cite{benmansour:2010dm, bonnivard:2014tw} where it was derived for the fast marching method.

\section{The Algorithm} \label{sec:algorithm}
In this section, we describe how to include the topological term in an explicit fashion in given finite element code. We note that the limiting effect on the time-step size of this explicit treatment is moderate, as the term does not depend on spatial gradients of the phase field function $u$.

The description is given in the two-dimensional case assuming that the finite element space corresponds to a triangulation of $\Omega$ with grid length scale $h$. Dimension three and more general basis element shapes can be treated by the same method. 

In the set up of the simulation, we create the dual graph $\Gamma_h$ corresponding to the finite element triangulation $\mathcal{T}_h$ as described in section~\ref{sec:disc_dist} (with the edge weights left unassigned for the time being, as they will change in each time step). Each triangle is furthermore associated with its volume $|T|$ and diameter $\diam(T)$.

Given a Galerkin space function $u = u^k$ in time step $k$, do the following.

\begin{enumerate}
\item\label{compute average} For all triangles $T$ in the triangulation, compute the average integral
\[
u_T = \frac1{|T|}\int_T u\dx.
\]

\item\label{compute weight} For each edge $e=e(T_1,T_2)$ in $\Gamma_h$ corresponding to two adjacent triangles $T_1$ and $T_2$, compute its weight as
\[
w_e = \frac{F(u_{T_1}) + F(u_{T_2})}{2}\,\cdot\,\frac{\diam(T_1) + \diam(T_2)}{2}.
\]
The second factor can be replaced by a generic grid length scale $h$ if the triangles are sufficiently uniform.

\item\label{find interesting elements} Create a list $I$ of all interface elements, i.e.\ all triangles such that
\[
u_T \in [\alpha,\beta].
\]

\item\label{create components} Separate the elements in $I$ into connected components $\{C_j\}_{j=1}^M$, where two triangles $T_1, T_2 \in I$ belong to the same component if $d^{\Gamma_h(u)}(T_1,T_2) = 0$. If there is only one connected component, the algorithm can be terminated here as our approximations of both $\C_\eps$ and its variation vanish.

\item\label{compute component integrals} For $j\in\{1,\dots,M\}$ calculate
\[
\overline{W}_j = \frac1\eps \sum_{T_l\in C_j} \widetilde W(u_{T_l})|T_l|.
\]

\item For $i,j \in \{1,\dots,M\}$, $i\neq j$ calculate the component distances $\overline{d}_{ij} = \dist^{\Gamma_h(u)}(C_i,C_j)$ as well as the shortest connecting paths between components $\overline{\gamma}_{ij} = \{T^{ij}_0, T^{ij}_1,\dots, T^{ij}_{L(i,j)}\}$, $T^{ij}_0 \in C_i$, $T^{ij}_{L(i,j)} \in C_j$. These computations can be performed using Dijkstra's algorithm on $\Gamma_h(u)$, see~\cite{dijkstra1959note}.

\item The approximate topological energy can now be computed as 
\[
\overline{\C}_\eps(u) = \sum_{i\neq j} \overline{d}_{ij}\overline{W}_i\overline{W}_j.
\]
We note that, compared to the original double integral term, this is a major simplification which is due to the specific choice of $F$ vanishing identically on the support of $\widetilde{W}$.

\item Our discrete approximation of the variation of $\overline{\C}_\eps$ with respect to a finite element basis function $\phi$ is then given by
\begin{align*}
&\delta_{u;\phi}\overline{\C_{\eps}}(u) = \sum_{i\neq j} 2\left[\frac{1}{\eps}\sum_{T_l\in C_i} \widetilde{W}'(u_{T_l})\int_{T_l}\phi\dx \right]\cdot \overline{W}_j\overline{d}_{ij}  
\\ &\quad+  \sum_{i\neq j} \overline{W}_i\overline{W}_j \cdot \delta_{u; \phi}\,\dist^{\Gamma_h}(C_i, C_j)
%\left[\sum_{\stackrel{T^{ij}_l \in \overline{\gamma}_{ij}}{l\neq 0}} \sum_{k=0}^1\frac{\partial}{\partial u_{T_{l-k}}}w_{e(T_{l-1},T_l)}(u_{T_{l-k}},u_{T_l}) \cdot \frac{1}{|T_{l-k}|}\int_{T_{l-k}}\phi\dx \right]
\end{align*}
where $\delta_{u; \phi}\,\dist^{\Gamma_h}(C_i, C_j) = \delta_{u; \phi}\,d^{\Gamma_h}(T_i, T_j)$ for any $T_i\in C_i$, $T_j\in C_j$ is given by \eqref{eq variation of distance}.
\end{enumerate}

This algorithm can be added to a given finite element implementation. We may compute the time-step from $u^{k}$ to $u^{k+1}$ with any scheme 
\[
\eps\langle u_{k+1} - u_k, \phi\rangle_{L^2} = \tau\left[ \Phi(u_{k+1}, u_k; \phi) - \frac{1}{\eps}\delta_{u;\phi}\overline{\C_{\eps}}(u) \right]
\]
that treats the topological term explicitly. Here $\tau$ is the time-step size and $\Phi$ is an explicit, implicit or mixed approximation of the variation $-\delta\F_\eps$.  In simulations, it has proven useful to use the sum of two functionals of type $\C_\eps$: one, to keep the portion of the interface close to $+1$ connected (i.e., with $\alpha$ and $\beta$ close to $+1$) and one to keep the portion of the interface close to $-1$ connected (i.e., with $\alpha$ and $\beta$ close to $-1$).

\section{Numerical Results} \label{sec:numerics}

We compare the discrete gradient flow~\eqref{eq:time_step} of $\F_\eps + \frac{a}{\eps}\C_\eps$, where $\F_\eps$ is given in~\eqref{eq:will_h0} and either $a=0$ or $a=6.0\cdot{10^1}$. The computational domain is a discretization of the three dimensional unit ball by approximately $1.6\cdot 10^{6}$ tetrahedral $P1$ finite elements. Our time-stepping algorithm is a simple first order fully implicit Euler scheme coupled to an explicit treatment of $\C_\eps$ as described in section~\ref{sec:algorithm}. The time step size is given by $\tau = 5\cdot 10^{-7}$ (with a smaller time step if $\tau = 1\cdot 10^{-8}$ for the first several hundred time steps during the fast convergence to the optimal profile). The higher spatial gradient in the energy is treated using a Ciarlet-Raviart-Monk mixed formulation~\cite{ciarlet1974mixed, Monk:1987ha} with clamped boundary conditions $u=-1$ on $\partial \Omega$ and $\frac{\partial u}{\partial n} = 0$ on $\partial \Omega$ and we note that on our convex domain this variational crime is only a misdemeanor~\cite{Gerasimov:2012fn}.

\begin{figure}[ht!]\begin{center}
\raisebox{-0.5\height}{\includegraphics[height=3.25cm]{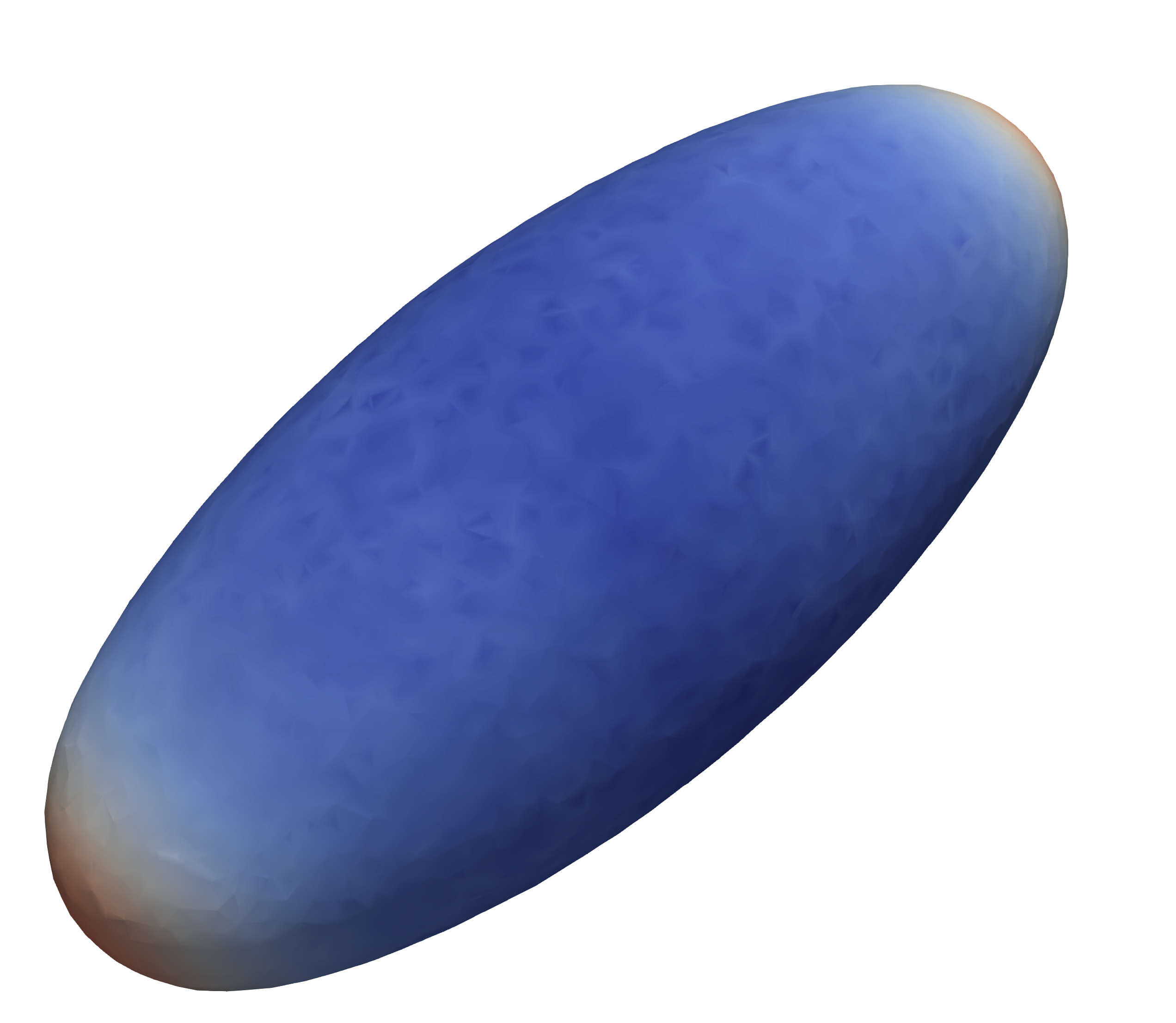}}
\raisebox{-0.5\height}{\includegraphics[height=3.25cm]{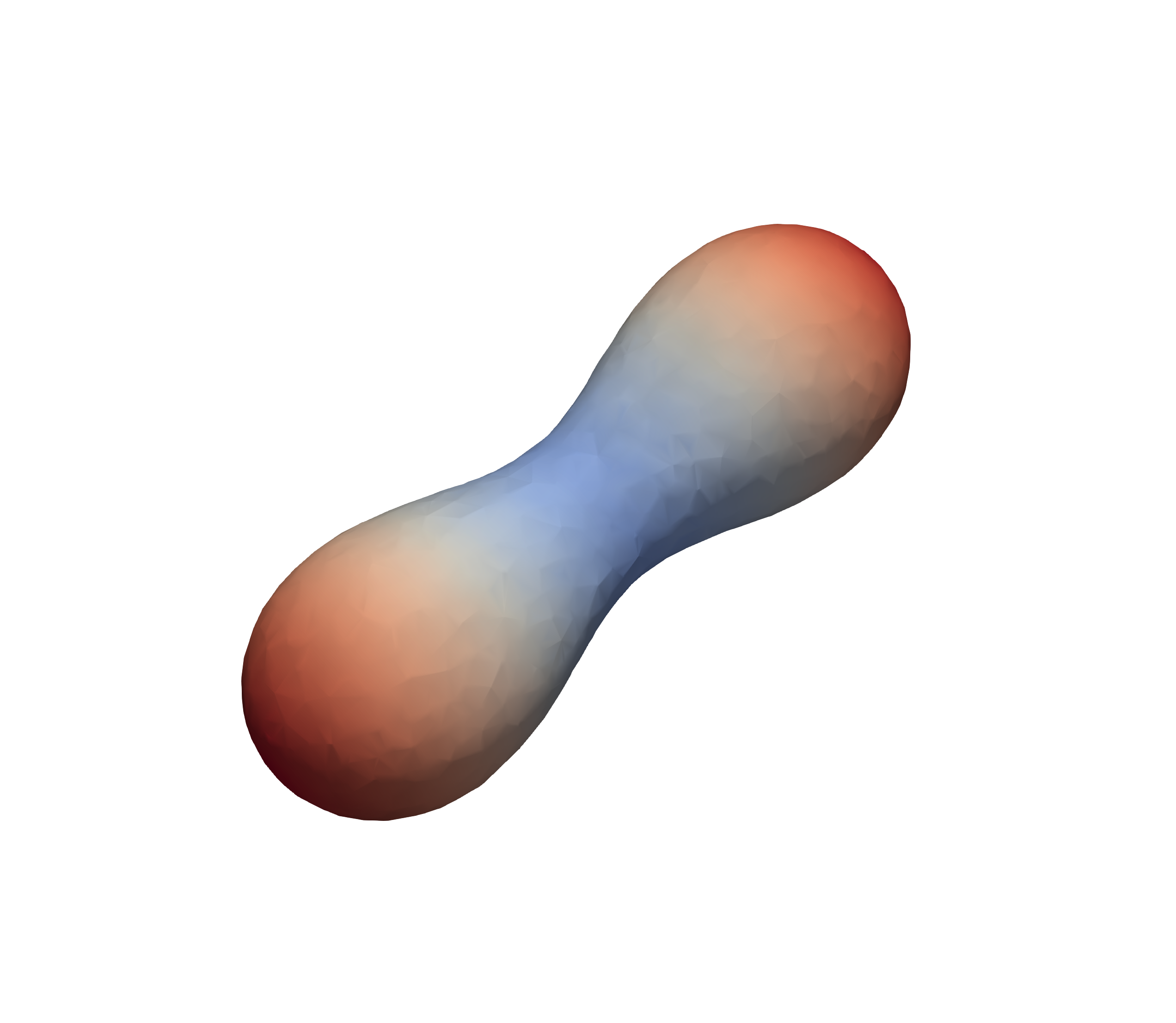}}
\raisebox{-0.5\height}{\includegraphics[height=3.25cm]{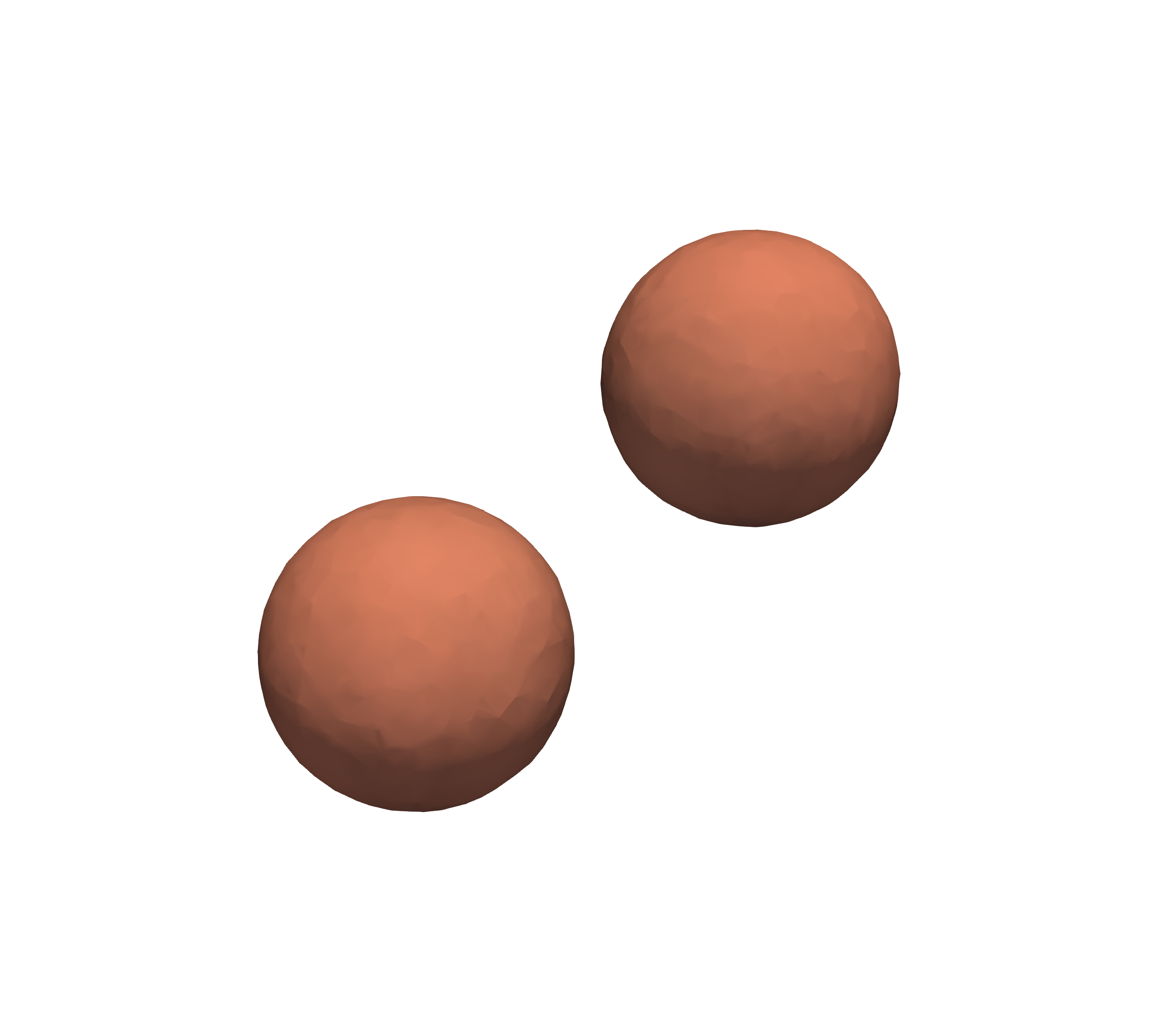}}
\raisebox{-0.5\height}{\includegraphics[height=2cm]{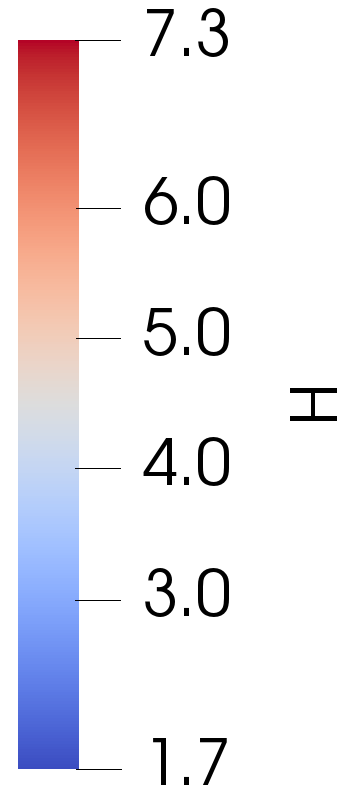}}
\caption{Simulation results for the geometric flow \emph{without} connectedness penalty. Shown are the initial condition (after some relaxation) at time $t=1\cdot 10^{-5}$, an intermediate configuration just before pinch-off at $t=3.0 \cdot 10^{-4}$ and the final state which was reached at $t=5.0\cdot 10^{-4}$. The color indicates the mean curvature, the image show the zero-level set of the phase field function $u$.}\label{fig:sim_will_nc}
\end{center}
\end{figure}
\begin{figure}[ht!]\begin{center}
\raisebox{-0.5\height}{\includegraphics[height=3.25cm]{init_0.png}}
\raisebox{-0.5\height}{\includegraphics[height=3.25cm]{mid_30.png}}
\raisebox{-0.5\height}{\includegraphics[height=3.25cm]{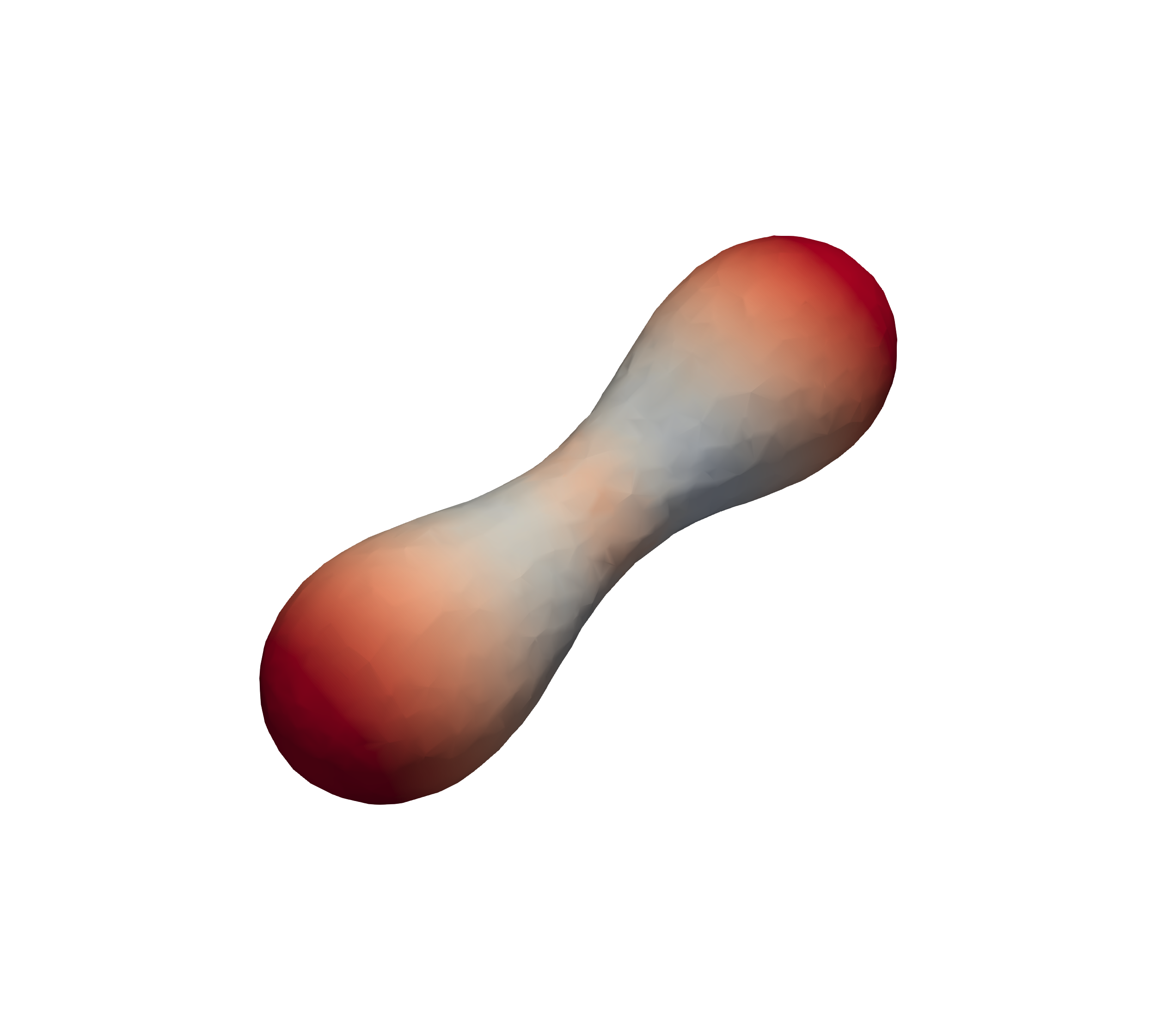}}
\raisebox{-0.5\height}{\includegraphics[height=2cm]{colorbar.png}}
\caption{Simulation results for the geometric flow \emph{with} connectedness penalty. The first two images are the same as in Figure~\ref{fig:sim_will_nc}, since no disconectedness had occurred yet. The final state in the third image was reached at $t=4.0\cdot 10^{-4}$.} \label{fig:sim_will_c}
\end{center}
\end{figure}

The parameters in the simulation are $\eps=0.03$, $\lambda = 0.1$, and $H_0 = 6$. The initial condition is an approximation of the characteristic function of an ellipsoid with principle axes $0.7$, $0.3$, and $0.3$. As seen in Figure~\ref{fig:sim_will_nc}, in the simulation of the case $a=0$, without topological constraint, the surface undergoes a pinch-off and the final steady state is given by two spheres of radius approximately equal $\frac{1}{6}$. A similar result of pinch-off was observed in~\cite{Du:2006hl}\footnote{The simulation in~\cite{Du:2006hl} was performed for $\lambda=0$. We note that our simulation produces virtually the same results for $\lambda=0$, but we present the case of positive $\lambda$ since none of the analytic results in~\cite{roger:2006ta} are admissible unless $S_\eps(u_\eps)$ remains uniformly bounded as $\eps\to 0$.}.

The simulation results for the case $a=6.0\cdot{10^1}$, i.e., including the topological constraint, are shown in Figure~\ref{fig:sim_will_c}. One can clearly see that the pinch-off into two components has been suppressed and a dumbbell-like shape is the final result. We conjecture that this shape is in fact a stationary point of the energy as no further motion was observed in the simulation even on a longer time-scale.

\begin{figure}[ht!]\begin{center}
\includegraphics[height=4cm]{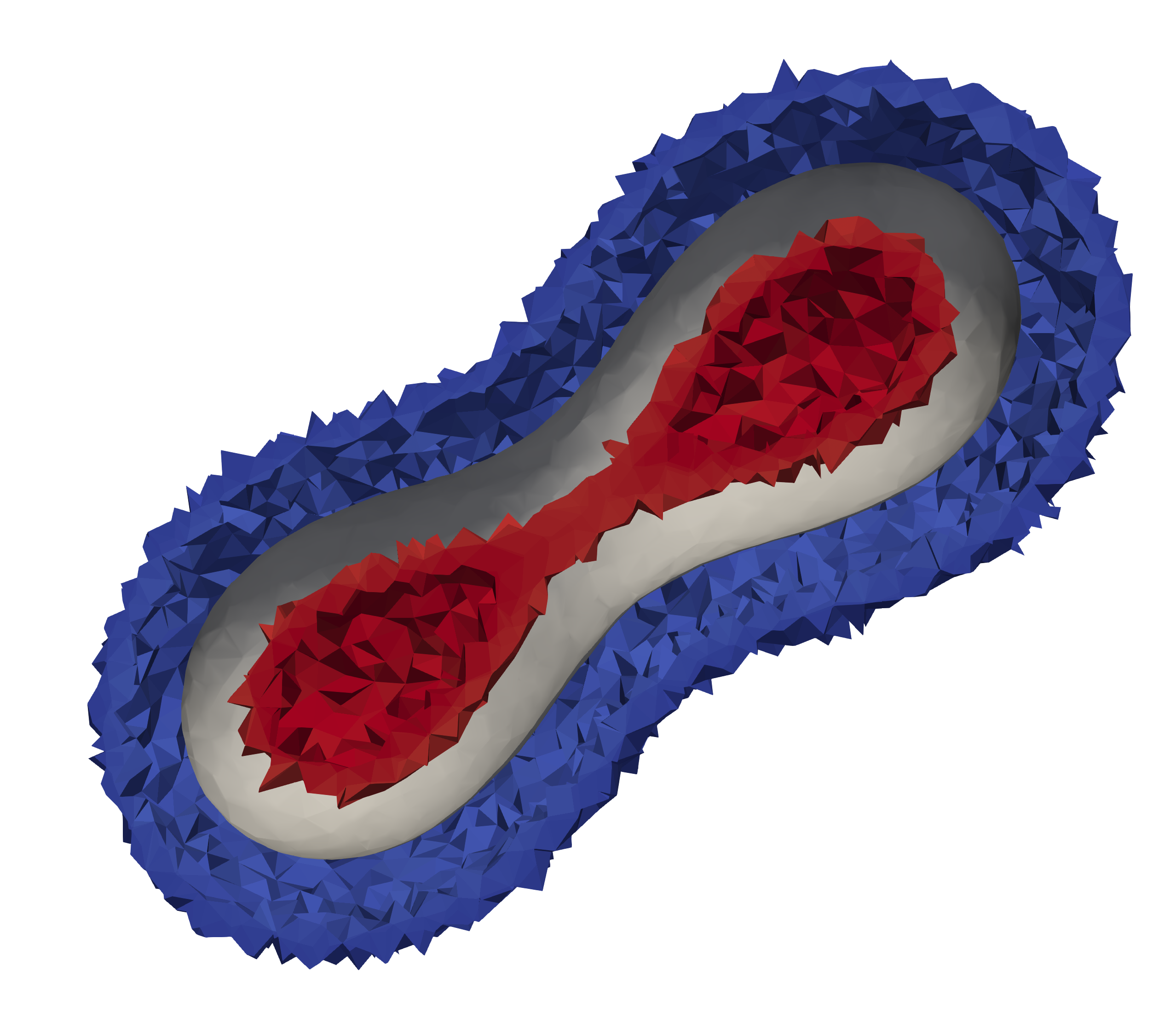}
\includegraphics[height=4cm]{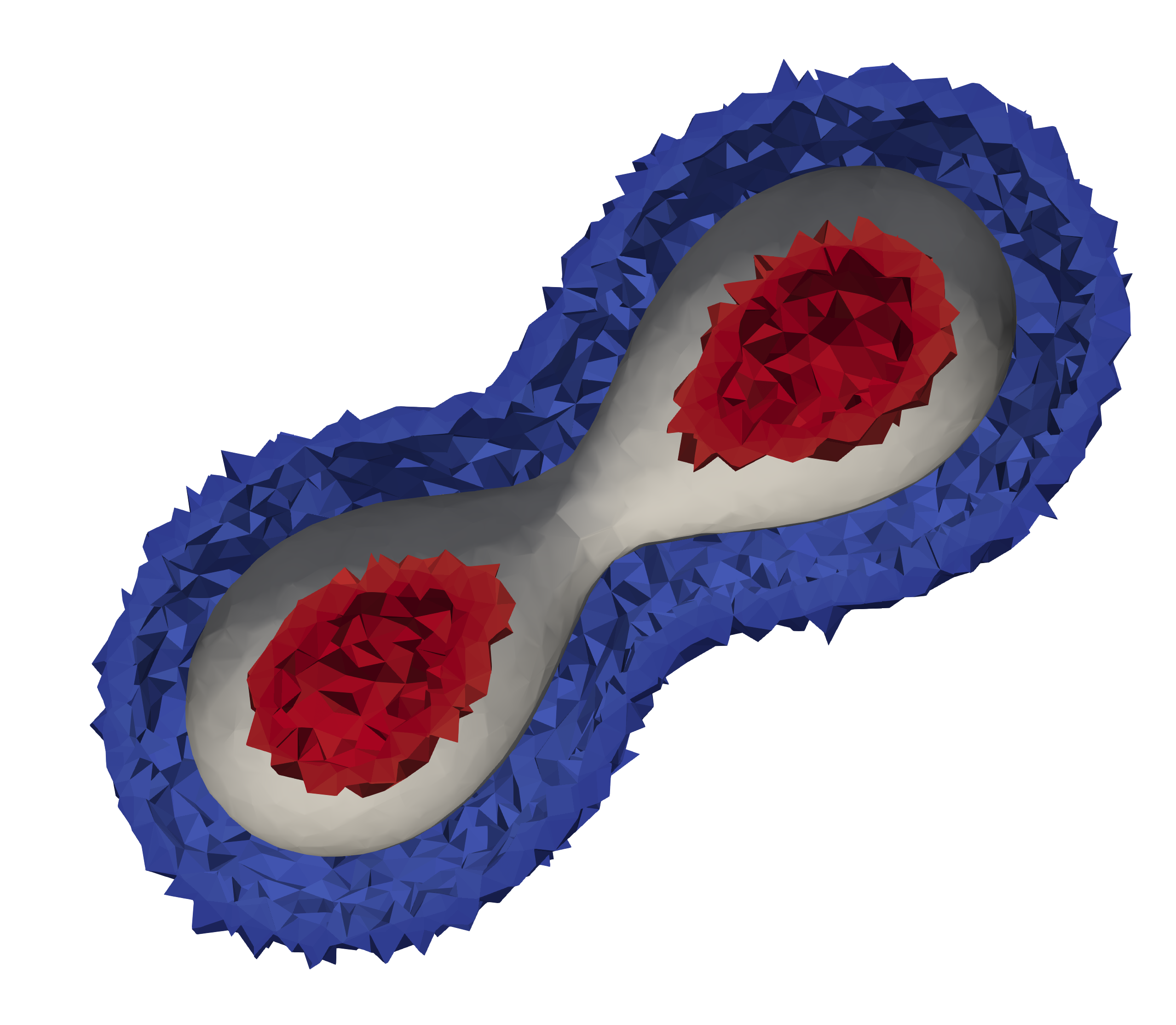}
\caption{The figures show the sets $\{-0.95 < u_T < -0.85\}$ (in blue) and $\{0.85 < u_T < 0.95\}$ (in red) as well as the zero-level set of $u$ (in grey). The image on the left is just before the start of disconnectedness is detected by the functional keeping the part of the transition layer close to $+1$ connected, at $t=3.0\cdot 10^{-4}$, the image on the left is shortly thereafter at $t=3.2\cdot 10^{-4}$. The functional $\C_\eps$ is zero in the left image, and non-zero in the right image.} \label{fig:shells}
\end{center}
\end{figure}

As mentioned at the end of section~\ref{sec:algorithm}, we implement the connectedness constraint using the sum of two functionals of type $\C_\eps$, in this specific case one with $\alpha=0.85$, $\beta=0.95$ (in order to keep the part of the transition layer close to the phase $u=+1$ connected) and another one with $\alpha=-0.95$, $\beta=-0.85$ (in order to keep the part of the transition layer close to the phase $u=-1$ connected). The functions $F$ and $\widetilde{W}$ in $\C_\eps$ are given by 
\begin{align}
F(s) &= \begin{cases}
(s-\alpha)^2\cdot c_1 & s<\alpha\\
0 & \alpha \le s \le \beta \\
(\beta-s)^2\cdot c_2 & s>\beta
\end{cases}
\quad\text{and}\\
\widetilde{W}(s) &= \begin{cases}
0 & s\le \alpha \\
(s-\alpha)^2(\beta-s)^2 \cdot c_3 & \alpha<s<\beta \\
0 & s \ge \beta
\end{cases},
\end{align}
 respectively, with $c_1$ and $c_2$ chosen such that $F(-1)=F(+1)=1$ and $c_3$ such that $\int_{-\infty}^\infty \widetilde{W}(s)\ds =1$. We note that in our case, at the pinch-off, the function $u$ dips below $0.85$ and thus the start of the interface becoming disconnected is detected by the algorithm. For an illustration see Figure~\ref{fig:shells}.

The equilibrium was reached in approximately 9 hours of wall-time using 8 cores of a computer server equipped with two Intel Xeon E5-2690 v4 processors. We note that the cpu-time spent computing the topological constraint is negligible (less than 0.1\% of the total cpu-time)---the only computational down-side may thus be the aforementioned restriction on the time-step size due to the necessary explicit treatment of $\C_\eps$. Implicit treatment of $\C_\eps$ is analytically questionable due to the lower regularity of $\C_\eps$.

\section{An application to image segmentation} \label{sec:2d_conn}

\begin{figure}[ht!]\begin{center}
\raisebox{-0.5\height}{\includegraphics[height=3.7cm]{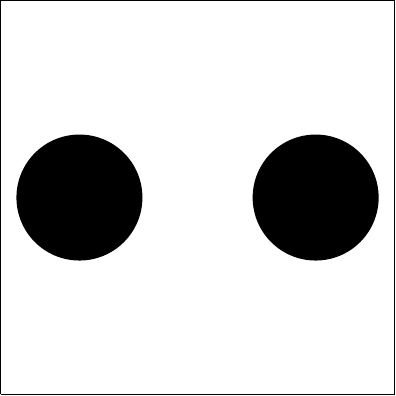}}
\raisebox{-0.5\height}{\includegraphics[height=3.7cm, clip, trim = 26.17cm 5.65cm 26.17cm 5.65cm]{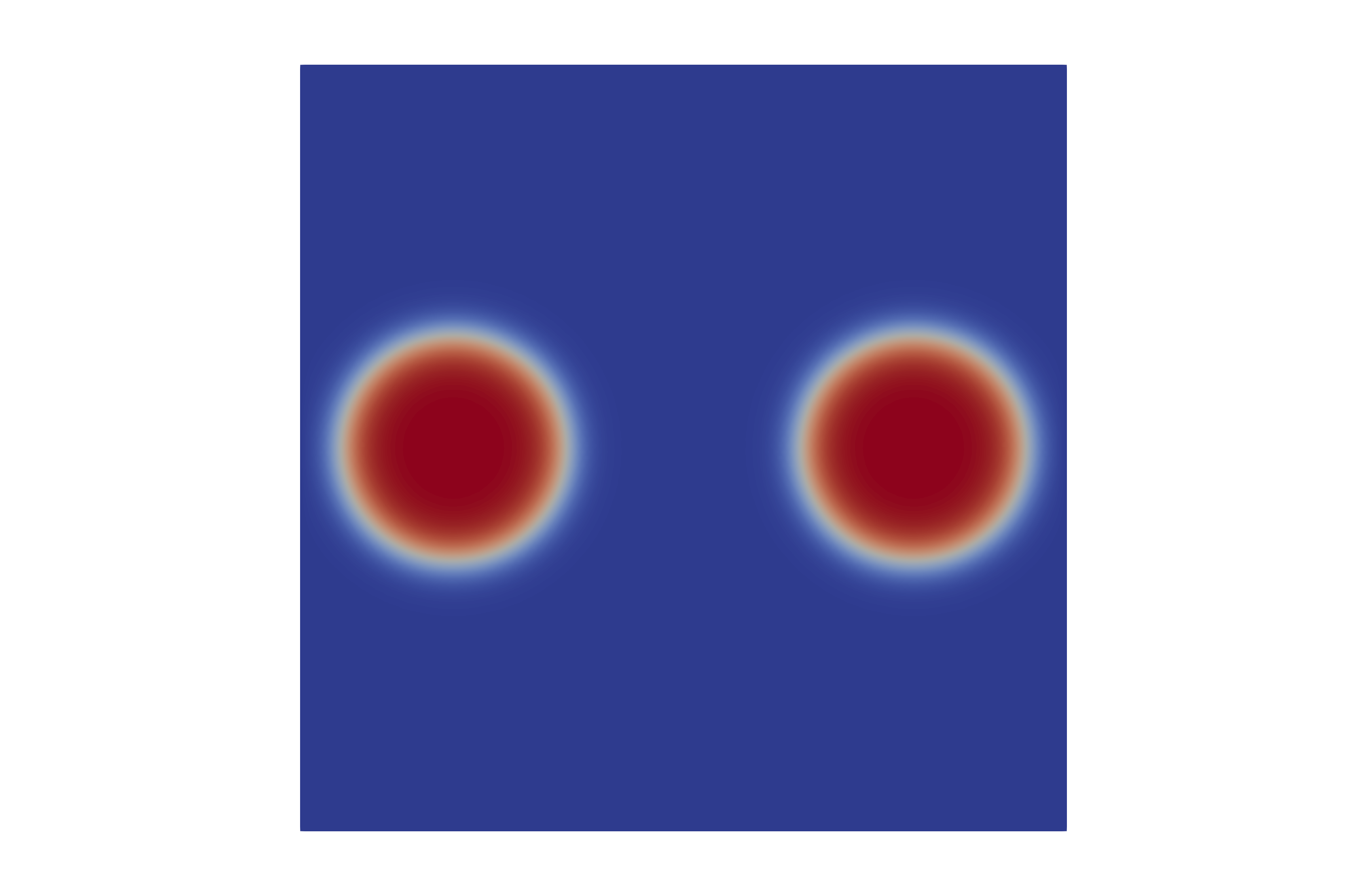}}
\raisebox{-0.5\height}{\includegraphics[height=3.7cm, clip, trim = 26.17cm 5.65cm 26.17cm 5.65cm]{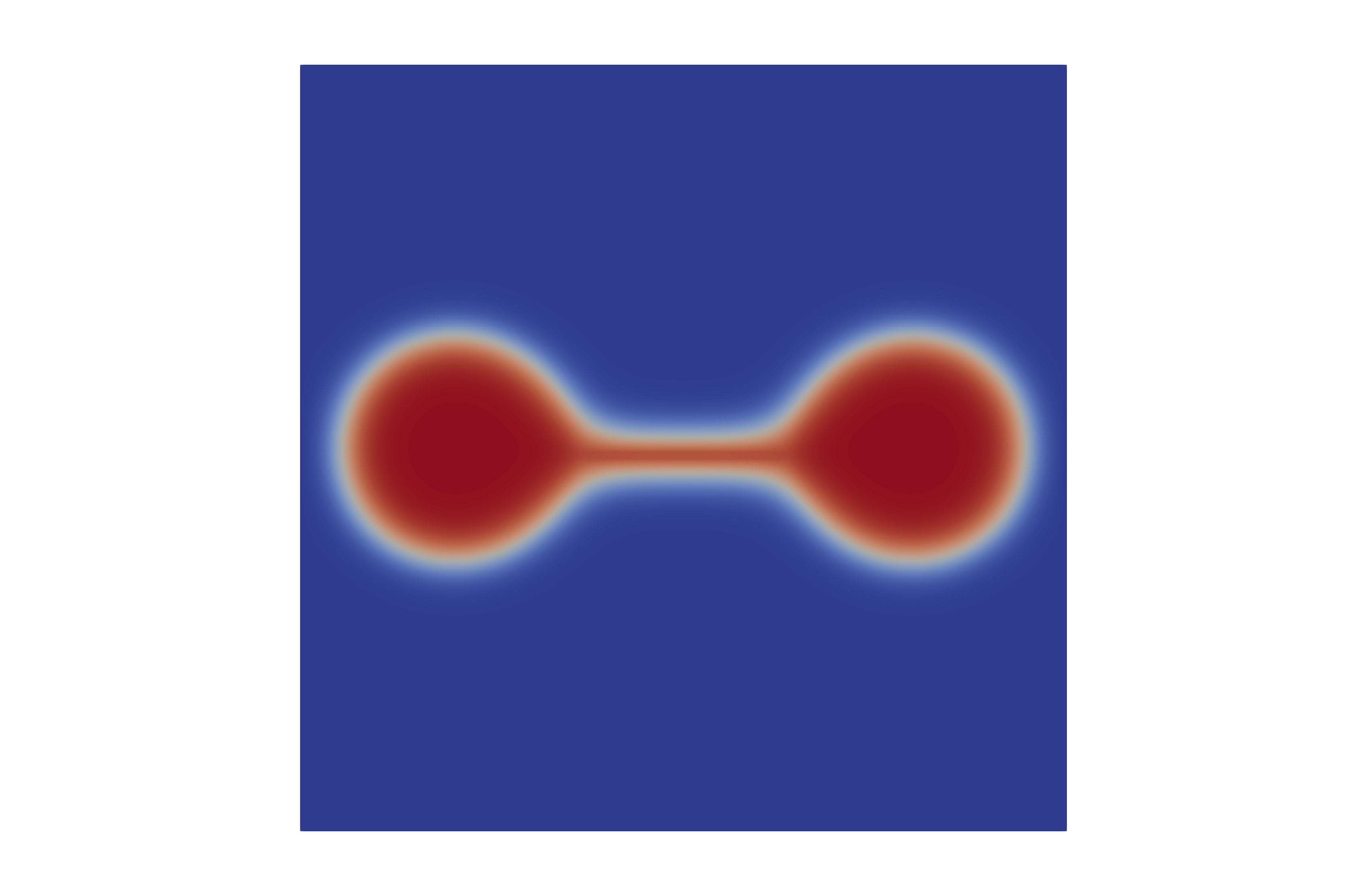}}
\raisebox{-0.5\height}{\includegraphics[height=2.5cm]{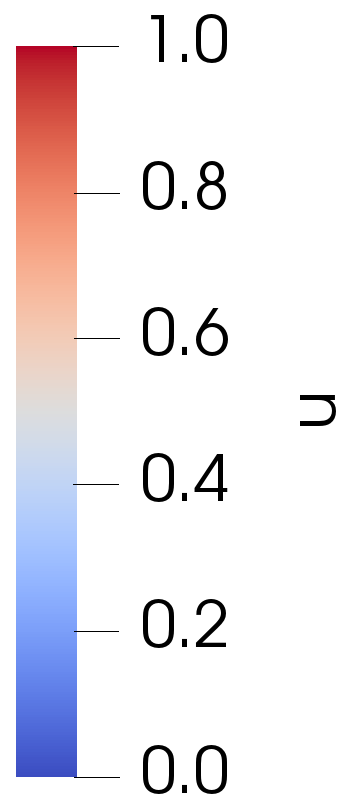}}
\caption{Stationary points for the image segmentation example with large disks (radius 0.16, distance of centers 0.6). From left to right: given black and white image $g$ (black corresponding to the value $+1$, in this case), stationary state $u$ without connectedness penalty, stationary state with connectedness penalty.}
\label{fig:image_seg1}
\end{center}
\end{figure}
\begin{figure}[ht!]\begin{center}
\raisebox{-0.5\height}{\includegraphics[height=3.7cm]{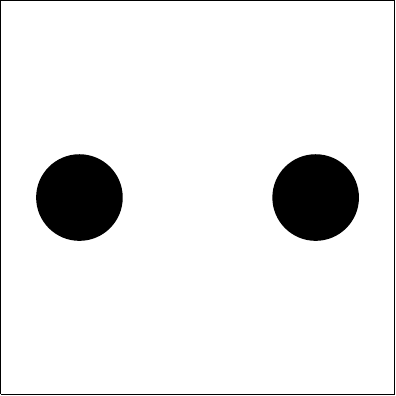}}
\raisebox{-0.5\height}{\includegraphics[height=3.7cm, clip, trim = 26.17cm 5.65cm 26.17cm 5.65cm]{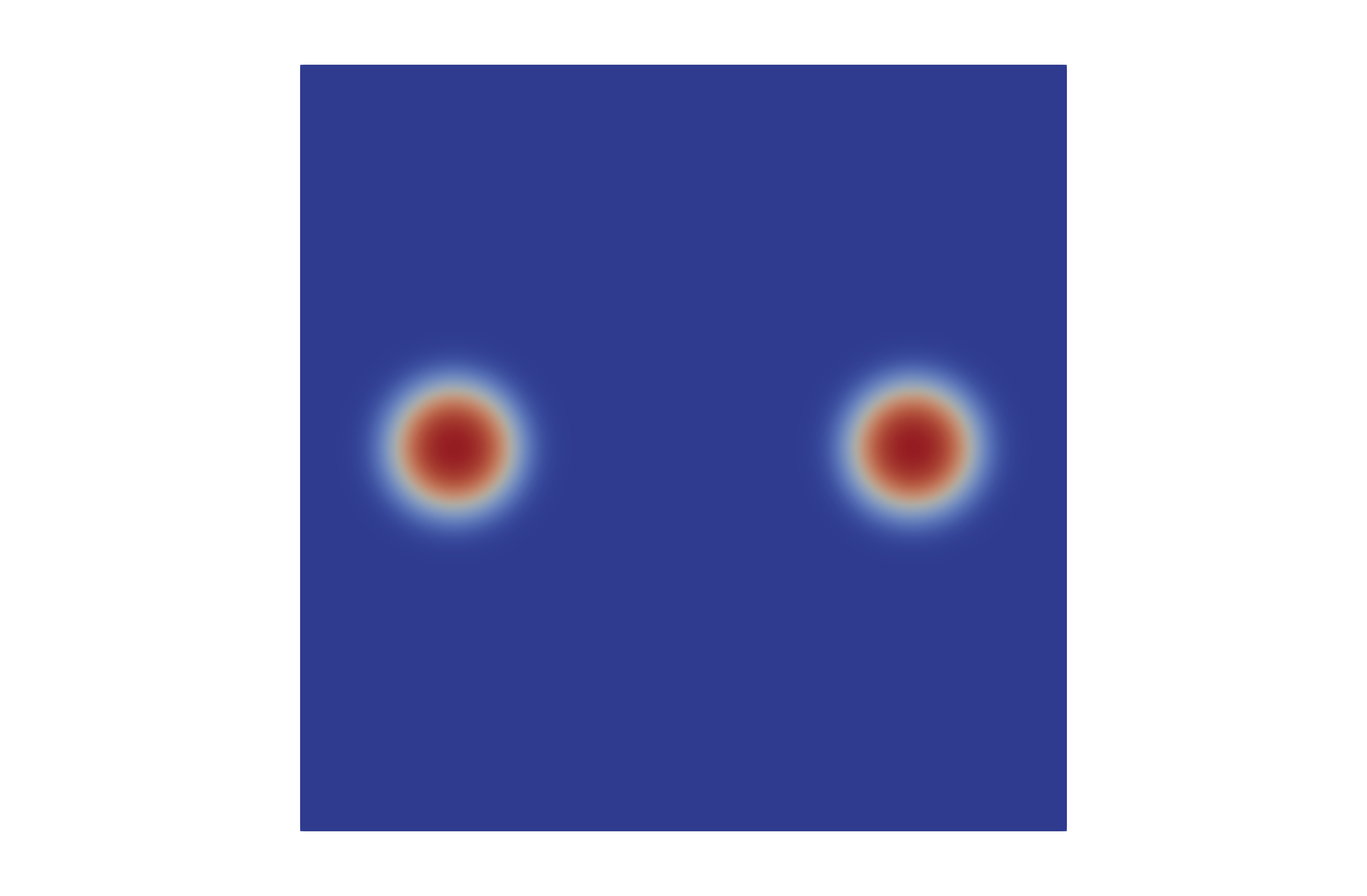}}
\raisebox{-0.5\height}{\includegraphics[height=3.7cm, clip, trim = 26.17cm 5.65cm 26.17cm 5.65cm]{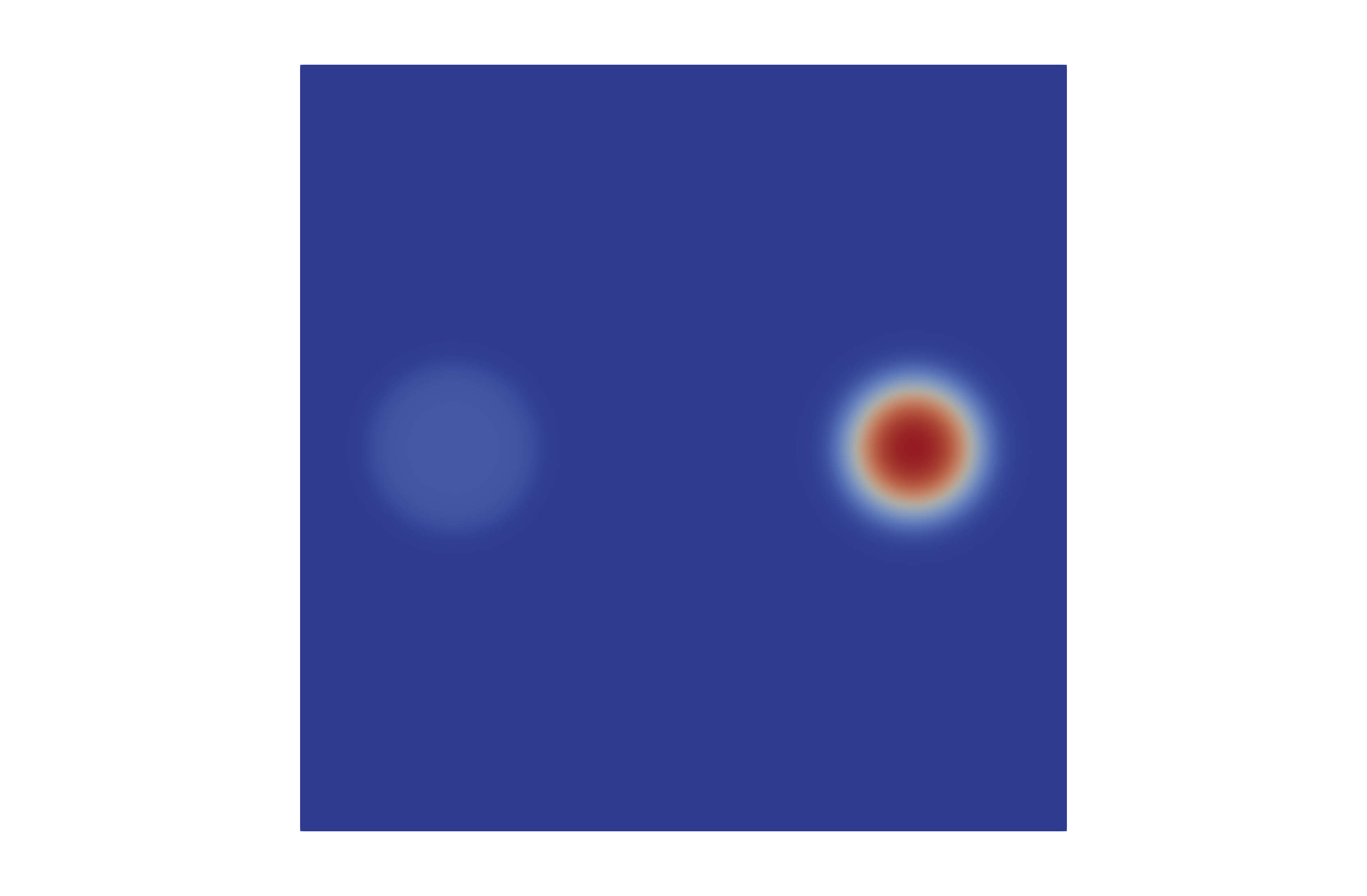}}
\raisebox{-0.5\height}{\includegraphics[height=2.5cm]{colorbar_u.png}}
\caption{Stationary points for the image segmentation example with small disks (radius 0.11, distance of centers 0.6). From left to right: given black and white image $g$, stationary state $u$ without connectedness penalty, stationary state with connectedness penalty.} \label{fig:image_seg2}
\end{center}
\end{figure}
To conclude, we give an outlook to a different application of the topological functional introduced in this article. We consider the energy
\[
\F_\eps(u) + \frac{a}{\eps^2} \,\C^\text{img}_\eps(u)= 
S_\eps(u)
+ \eta\int_\Omega |u-g|^2\dx+ \frac{a}{\eps^2}\,\C^\text{img}_\eps(u),
\]
with $\eta>0$ and $g\colon \Omega \to [0,1]$ given and $\Omega \subset \R^2$. Again, the case $a=0$ is that of no disconnectedness penalty, and for $a>0$ disconnectedness of the set $\{\alpha < u < \beta\}$ is penalized.

For $a=0$ and $S_\eps$ using a usual double-well energy $W(u)=\frac{1}{4}u^2(u-1)^2$ (in this case with minima at $0$ and $1$ and normalization constant $c_0 = \frac{\sqrt{2}}{12}$), the functional is a typical image-segmentation functional with perimeter regularization and a fidelity term.

For $a>0$ we now choose $a=4.0\cdot 10^{-1}$ and $\alpha = 0.9$ and $\beta=1.2$\footnote{The functions $\widetilde{W}$ and $F$, as well as the constants $c_1$ and $c_3$ are picked as before, the constant $c_2$ is irrelevant, noting that in the simulation it always holds that $u<1$. The functional, $\C^\text{img}_\eps$ is given by $\eps^2\C_\eps$, since the integrals are now performed over open sets and not over boundary layers.}: this means that the topological term will enforce path-connectedness of the set $\{u\approx +1\}$ and so the grey scale image given by $g$ should be segmented into a \emph{connected} segment and its exterior.

In the present example, $\Omega=\left(-\frac{1}{2},\frac{1}{2}\right)^2$ is the unit square, $\eps=1\cdot 10^{-2}$ and $\eta = 1.05\cdot 10^1$. The discretization is given by approximately $2.3\cdot 10^{4}$ $P1$ triangular finite elements. We again compute a time-discrete $L^2$-gradient flow (using now a semi-implicit first order Euler scheme with only the linear highest gradient term being treated implicitly) of the energy until a stationary state has been reached. In all cases, the initial condition is given by the characteristic function of the set $\{r < 0.25+0.15\cos(5.0\,\theta)\}$ in polar coordinates.

In Figures~\ref{fig:image_seg1} and~\ref{fig:image_seg2} the results are shown. On the left, the source image to be segmented is displayed (in our simple example, $g$ itself only takes values in $\{0,1\}$ and is the characteristic function of two separated disks). The two other pictures show the computed phase field minimizer $u$, first without connectedness penalty ($a=0$) then with connectedness penalty ($a>0$). One can clearly see that without this penalty, the source image is simply reproduced.

\begin{table}[ht!]
\begin{center}
\begin{tabular}{l|l|l}
 & large disks & small disks \\
\hline
Radii &0.16 & 0.11\\
\hline
Distance between centers of disks & 0.6 & 0.6 \\
\hline
Twice the distance between disks & 0.56 & 0.76\\
\hline
Fidelity penalty for removing one disk & 0.84 & 0.40
\end{tabular}
\caption{Energy comparison for the cost of connectedness (either by adding a double layer between connected components of the image or by removing a connected component).} \label{tab:energies}
\end{center}
\end{table}

In Figure~\ref{fig:image_seg1} connectedness is established by adding a double-layer between the two disks in the image $g$. With the disks being smaller in Figure~\ref{fig:image_seg2} and our choice of parameters, the double layer would be more energetically costly than the error being made in the fidelity term, so one of the disks is simply being ignored. The respective values for the energies have been assembled in Table~\ref{tab:energies}. We do note, however, that a global optimum cannot generally be found by such a gradient flow as local minimal of the energy exist (in our setting, the global minimizers are given by $u\approx 0$ everywhere). A rigorous analysis of this image-segmentation energy with connectedness constraint will be presented in~\cite{novaga-etus}.

\section*{Acknowledgements} PWD gratefully acknowledges partial support by the Wissenschaftler-R{\"u}ckkehr\-programm GSO/CZS as well as inspiring discussions with Benedikt Wirth (M{\"u}nster), Sebastian Reuther (Dresden), and Douglas N.~Arnold (Minneapolis).

%\bibliographystyle{alphaabbr}
%\bibliography{impl}

\end{document}